\documentclass[11pt, reqno]{amsart}
\usepackage{amssymb}
\usepackage{amsmath}
\usepackage{amsthm}
\usepackage{amsfonts}
\usepackage[latin5]{inputenc}

\newtheorem{theorem}{Theorem}
\newtheorem{corollary}{Corollary}
\newtheorem{lemma}{Lemma}

\begin{document}
\title[Rotational hypersurfaces family satisfying $\mathbb{L}_{n-3}\mathcal{G%
}=\mathcal{A}\mathcal{G}$ in ${\mathbb{E}}^{n}$]{Rotational hypersurfaces
family satisfying $\mathbb{L}_{n-3}\mathcal{G}=\mathcal{A}\mathcal{G}$ in
the $n$-dimensional Euclidean space}
\author[E. G\"{u}ler]{Erhan G\"{u}ler$^{\ast }$}
\address[Erhan G\"{u}ler]{ Bart\i n University, Faculty of Sciences,
Department of Mathematics, Kutlubey Campus, 74100 Bart\i n, Turkey}
\email{eguler@bartin.edu.tr}
\author[N. C. Turgay]{Nurettin Cenk Turgay}
\address[Nurettin Cenk Turgay]{ Istanbul Technical University, Faculty of
Science and Letters, Department of Mathematics, 34469 Maslak, Istanbul,
Turkey }
\email{turgayn@itu.edu.tr}
\thanks{$^{\ast }$Corresponding Author}
\date{}
\dedicatory{}
\thanks{This paper is in final form and no version of it will be submitted
for publication elsewhere}
\date{\today }

\begin{abstract}
In this paper, we investigate rotational hypersurfaces family in $n$%
-dimensional Euclidean space ${\mathbb{E}}^{n}$. Our focus is on studying
the Gauss map $\mathcal{G}$ of this family with respect to the operator $%
\mathbb{L}_{k}$, which acts on functions defined on the hypersurfaces. The
operator $\mathbb{L}_{k}$ can be viewed as a modified Laplacian and is known
by various names, including the Cheng--Yau operator in certain cases.
Specifically, we focus on the scenario where $k=n-3$ and $n\geq 3$. By
applying the operator $\mathbb{L}_{n-3}$ to the Gauss map $\mathcal{G}$, we
establish a classification theorem. This theorem establishes a connection
between the $n\times n$ matrix $\mathcal{A}$, and the Gauss map $\mathcal{G}$
through the equation $\mathbb{L}_{n-3}\mathcal{G}=\mathcal{A}\mathcal{G}$.
\end{abstract}

\maketitle

\section{Introduction}

Chen \cite{C} introduced the challenge of categorizing finite type surfaces
within the three dimensional Euclidean space $\mathbb{E}^{3}$. A Euclidean
submanifold is designated as having Chen's finite type characteristic when
its coordinate functions can be expressed as a limited combination of its $%
\Delta $ Laplacian (or Laplace--Beltrami operator)'s eigenfunctions \cite%
{C00, C}. Additionally, the concept of finite type can be expanded to
encompass smooth functions on a submanifold existing within either a
Euclidean space or a pseudo-Euclidean space. Readers can refer to \cite%
{Chen1991, C2, CI1993, CP1987, DPV, HV1992, M, M2, T, Y} for other studies
related to the subject mentioned above.

On the other hand, the extended version of the Laplace--Beltrami operator $%
\mathbb{L}_{0}=\Delta $ is referred to as the Cheng--Yau operator $\mathbb{L}%
_{1}=\square $. The general operator is denoted by the $\mathbb{L}_{k}$
operator.

Numerous investigations carried out concerning the previously mentioned
topics. For instance, Kim and Turgay \cite{KimTurgay} studied surfaces with $%
\mathbb{L}_{1}$-pointwise $1$-type Gauss map in ${\mathbb{E}}^{4}$.
Kim\textperiodcentered et al. \cite{Kim et al} focused on $\mathbb{L}_{1}$
operator and Gauss map of surfaces of revolution in $\mathbb{E}^{3}$. G\"{u}%
ler and Turgay \cite{GulerTurgay} studied $\mathbb{L}_{1}$ operator and
Gauss map of rotational hypersurfaces in $\mathbb{E}^{4}.$ Mohammadpouri,
Kashani, and Pashaie \cite{MKP2013} introduced $\mathbb{L}_{1}$-finite type
Euclidean surfaces. Kashani \cite{Kashani2009} studied $\mathbb{L}_{1}$%
-finite type (hyper)surfaces in $\mathbb{E}^{n+1}.$

{Alias and }Kashani \cite{AliasKashani2010} gave hypersurfaces in space
forms satisfying the condition $\mathbb{L}_{k}x=Ax+b,$ where $A\in \mathbb{R}%
^{\left( n+2\right) \times \left( n+2\right) }$ is a constant matrix$,$ $%
b\in \mathbb{R}^{n+2}$ is a constant vector. Lucas and Ramirez-Ospina \cite%
{LR2011, LR2012} worked hypersurfaces in the Lorentz-Minkowski space
satisfying $\mathbb{L}_{k}\psi =A\psi +b.$  See also 
\cite{LR2013a, LR2013b} for their works. Mohammadpouri and Kashani \cite%
{MK2013} studied quadric hypersurfaces of $\mathbb{L}_{r}$-finite type.
Pashaie and Kashani \cite{PashaieKashani2013} considered spacelike
hypersurfaces in Riemannian or Lorentzian space forms, and they also\cite{PashaieKashani2014} gave timelike hypersurfaces in the Lorentzian space
forms, satisfying $\mathbb{L}_{k}x=Ax+b.$ Mohammadpouri \cite%
{Mohammadpouri2018} introduced the hypersurfaces with $\mathbb{L}_{r}$-pointwise 1-type Gauss map.

The aim of this paper is to investigate the properties of a family of
rotational hypersurfaces in $\mathbb{E}^{n}$ and establish a connection
between the Gauss map $\mathcal{G}$ and specific geometric objects, such as
hyperplanes, right circular hypercones, circular hypercylinders, and
hyperspheres, by analyzing the equation $\mathbb{L}_{k}\mathcal{G}=\mathcal{%
AG}$ for a $n\times n$ matrix $\mathcal{A}$, $k=n-3\geq
0$, with integers $n$, involving $\mathbb{L}_{k}$ operators.

For further exploration of the $\mathbb{L}_{k}$ operators mentioned
previously, readers are encouraged to refer to the comprehensive
investigation conducted by Chen et al. \cite{CGYH}.

We present the main theorem of this work in this section. We then begin by
recalling the fundamental notions of $n$-dimensional Euclidean geometry in
Section 2. Moving on to Section 3, we review the $\mathbb{L}_{k}$ operators.
Then, in Section 4, we provide the definition of a rotational hypersurfaces
family in Euclidean $n$-spaces. In Section 5, we calculate the relationships
between the mean curvature and the Gauss--Kronecker curvature of the
aforementioned family. The pivotal part of our study lies in Section 6,
where we apply the operator $\mathbb{L}_{n-3}$ to the Gauss map $\mathcal{G}$%
, leading us to establish a classification theorem. This theorem establishes
a connection between the $n\times n$ matrix $\mathcal{A}$ and the Gauss map $%
\mathcal{G}$ through the equation $\mathbb{L}_{n-3}\mathcal{G}=\mathcal{A}%
\mathcal{G}$ for the family with integers $n\geq 3$. 

Throughout this research, our main focus revolve around the following
theorem, which we prove in Section 6.

\textbf{Main} \textbf{Theorem.} \textit{Let }$\mathfrak{x}=\mathfrak{x}%
(r,\theta _{1},\theta _{2},\hdots,\theta _{n-2})$\textit{\ be a family of
rotational hypersurfaces} \textit{in} $\mathbb{E}^{n}$ \textit{given by}%
\begin{eqnarray*}
\mathfrak{x} &=&\left( f(r)\prod\limits_{i=1}^{n-2}\cos \theta
_{i},f(r)\sin \theta _{1}\prod\limits_{i=2}^{n-2}\cos \theta _{i},\right. \\
&&\text{ \ \ \ \ \ \ \ \ \ \ \ }\left. f(r)\sin \theta
_{2}\prod\limits_{i=3}^{n-2}\cos \theta _{i},\hdots,f(r)\sin \theta
_{n-3}\cos \theta _{n-2},f(r)\sin \theta _{n-2},\varphi (r)\right) .
\end{eqnarray*}%
\textit{Then, the Gauss map }$\mathcal{G}$ \textit{associated with} $%
\mathfrak{x}$\textit{\ satisfies the equation }$\mathbb{L}_{n-3}\mathcal{G}=%
\mathcal{AG}$ \textit{for a} $n\times n$ \textit{matrix} $\mathcal{A}$ 
\textit{with integers} $n\geq 3$\textit{,} \textit{if and only if }$%
\mathfrak{x}$ \textit{is an open part of the following}:

$(1)$ \textit{a hyperplane,}

$(2)$\textit{\ a right circular hypercone,}

$(3)$\textit{\ circular hypercylinder,}

$(4)$\textit{\ a hypersphere.}

\section{Preliminaries}

In this section, we describe notations, basic facts and definitions. Let $%
\mathbb{E}^{n}$ denote the Euclidean $n$-space with the canonical Euclidean
metric tensor given by $\widetilde{g}=\langle \ ,\ \rangle
=\sum\limits_{i=1}^{n}dx_{i}^{2},$ where $(x_{1},x_{2},\hdots,x_{n})$ is a
rectangular coordinate system in $\mathbb{E}^{n}$.

Consider an $n$-dimensional Riemannian submanifold of the space $\mathbb{E}%
^{n}$. We denote Levi--Civita connections of $\mathbb{E}^{n}$ and $M$ by $%
\widetilde{\nabla }$ and $\nabla $, respectively. We shall use letters $%
X,\;Y,\;Z,\;W$ (resp., $\xi ,\;\eta $) to denote vectors fields tangent
(resp., normal) to $M$. The Gauss and Weingarten formulas are given,
respectively, by 
\begin{eqnarray}
\widetilde{\nabla }_{X}Y &=&\nabla _{X}Y+h(X,Y),  \label{MEtomWeingarten} \\
\widetilde{\nabla }_{X}\xi &=&-\mathbf{A}_{\xi }(X)+D_{X}\xi ,
\end{eqnarray}%
where $h$, $D$ and $\mathbf{A}$ are the second fundamental form, the normal
connection and the shape operator of $M$, respectively.

For each $\xi \in T_{p}^{\bot }M$, the shape operator $\mathbf{A}_{\xi }$ is
a symmetric endomorphism of the tangent space $T_{p}M$ at $p\in M$. The
shape operator and the second fundamental form are related by 
\begin{equation*}
\left\langle h(X,Y),\xi \right\rangle =\left\langle \mathbf{A}_{\xi
}X,Y\right\rangle .
\end{equation*}

The Gauss and Codazzi equations are given, respectively, by 
\begin{eqnarray}
\langle R(X,Y,)Z,W\rangle &=&\langle h(Y,Z),h(X,W)\rangle -\langle
h(X,Z),h(Y,W)\rangle ,  \label{MinkCodazzi} \\
(\bar{\nabla}_{X}h)(Y,Z) &=&(\bar{\nabla}_{Y}h)(X,Z),
\end{eqnarray}%
where $R,\;R^{D}$ are the curvature tensors associated with connections $%
\nabla $ and $D$, respectively, and $\bar{\nabla}h$ is defined by 
\begin{equation*}
(\bar{\nabla}_{X}h)(Y,Z)=D_{X}h(Y,Z)-h(\nabla _{X}Y,Z)-h(Y,\nabla _{X}Z).
\end{equation*}

Now, let $M$ be an oriented hypersurface in the Euclidean space $\mathbb{E}%
^{n}$, $\mathbf{A}$ its shape operator and $x$ its position vector. We
consider a local orthonormal frame field $\{e_{1},e_{2},\hdots,e_{n}\}$ of
consisting of principal direction of $M$ corresponding from the principal
curvature $\kappa _{i}$ for $i=1,2,\hdots n$. Let the dual basis of this
frame field be $\{\theta _{1},\theta _{2},\hdots,\theta _{n}\}$. Then, the
first structural equation of Cartan is indicated by 
\begin{equation}
d\theta _{i}=\sum\limits_{i=1}^{n}\theta _{j}\wedge \omega _{ij},\quad i=1,2,%
\hdots,n  \label{CartansFirstStructural}
\end{equation}%
where $\omega _{ij}$ denotes the connection forms corresponding to the
chosen frame field. We denote the Levi-Civita connection of $M$ and $\mathbb{%
E}^{n+1}$ by $\nabla $ and $\widetilde{\nabla },$ respectively. Then, from
the Codazzi equation \eqref{MinkCodazzi}, the following hold 
\begin{eqnarray}
e_{i}(\kappa _{j}) &=&\omega _{ij}(e_{j})(\kappa _{i}-\kappa _{j}),
\label{Coda2} \\
\omega _{ij}(e_{l})(\kappa _{i}-\kappa _{j}) &=&\omega _{il}(e_{j})(\kappa
_{i}-\kappa _{l})
\end{eqnarray}%
for distinct $i,j,l=1,2,\hdots,n$. Considering 
\begin{equation*}
s_{j}=\sigma _{j}(\kappa _{1},\kappa _{2},\hdots,\kappa
_{n})=\sum\limits_{1\leq i_{1}<i_{2}<\hdots<i_{j}\leq n}\kappa
_{i_{1}}\kappa _{i_{2}}\hdots\kappa _{i_{j}},
\end{equation*}%
we use the notation 
\begin{equation*}
r_{i}^{j}=\sigma _{j}(\kappa _{1},\kappa _{2},\hdots,\kappa _{i-1},\kappa
_{i+1},\kappa _{i+2},\hdots,\kappa _{n}).
\end{equation*}%
By the definition, we have $r_{i}^{0}=1$ and $s_{n+1}=s_{n+2}=\cdots =0$.

On the other hand, we will call the function $s_{k}$ as the $k$-th mean
curvature of $M$. We would like to note that functions $H=\frac{1}{n}s_{1}$
and $K=s_{n}$ are called the mean curvature and Gauss--Kronecker curvature
of $M$, respectively. In particular, $M$ is said to be $j$-minimal if $%
s_{j}\equiv 0$ on $M$.

\section{$\mathbb{L}_{k}$ operators}

Let $M$ be a hypersurface of $\mathbb{E}^{n+1}$ with $\kappa _{1}(x),\hdots%
,\kappa _{n}(x)$ as its principal curvatures. Associated with the principal
curvatures, there are $n$ algebraic invariants given by 
\begin{equation*}
s_{k}(x)=\sigma _{k}(\kappa _{1}(x),\hdots,\kappa _{n}(x)),\quad 1\leq k\leq
n,
\end{equation*}%
where $\sigma _{k}:\mathbb{R}^{n}\rightarrow \mathbb{R}$ is the elementary
symmetric function defined by 
\begin{equation*}
\sigma _{k}(x_{1},\hdots,x_{n})=\sum_{i_{1}<\cdots <i_{r}}x_{i_{1}}\hdots %
x_{i_{k}}.
\end{equation*}%
Then the characteristic polynomial of the shape operator \textbf{$A$} of $M$
can be expressed in terms of the $s_{k}$'s as 
\begin{equation}
Q_{\mathbf{A}}(t)=\det (t\mathcal{I}-\mathbf{A})=%
\sum_{k=0}^{n}(-1)^{k}s_{k}t^{n-k},\;\;s_{0}=1,  \label{poly}
\end{equation}%
where $\mathcal{I}$ describes the identity map. The $k$\textit{-th mean
curvature} $H_{k}$ of $M$ is $s_{k}={\binom{{n}}{{k}}}H_{k},\,0\leq k\leq n.$
In fact, $H_{1}$ is the mean curvature, $H_{k}$ is intrinsic for even $k$,
and $H_{k}$ is extrinsic for odd $k$.

The Newton transformations $\mathcal{P}_{k}:\mathfrak{X}(M)\rightarrow 
\mathfrak{X}(M)$ $(k=1,\hdots,n)$ are defined inductively from the shape
operator $\mathbf{A}$ by 
\begin{equation}
\mathcal{P}_{0}=\mathcal{I},\hdots,\mathcal{P}_{k}=s_{k}\mathcal{I}-\mathbf{A%
}\circ \mathcal{P}_{k-1}={\binom{{n}}{{k}}}H_{k}\mathcal{I}-\mathbf{A}\circ 
\mathcal{P}_{k-1},  \label{Pk}
\end{equation}%
where $\mathcal{I}$ denotes the identity map on $\mathfrak{X}(M)$. By
applying Cayley--Hamilton's theorem, we have $\mathcal{P}_{n}=0$ from %
$\left(\ref{poly}\right)$. If $k$ is even, $\mathcal{P}_{k}$ does not depend on the chosen
orientation, but if $k$ is odd there is a change of sign in $\mathcal{P}_{k}$%
.

Associated with each Newton transformation $\mathcal{P}_{k}$, one has the
linear differential operator $\mathbb{L}_{k}:\mathcal{F}(M)\rightarrow {%
\mathcal{F}}(M)$ given by $\mathbb{L}_{k}(f)=-\mathrm{Tr}(\mathcal{P}%
_{k}\circ \nabla ^{2}f),$ where $\nabla ^{2}f:\mathfrak{X}(M)\rightarrow 
\mathfrak{X}(M)$ is the self-adjoint linear operator metrically equivalent
to the Hessian of $f$ and given by 
\begin{equation*}
\left\langle \nabla ^{2}f(X),Y\right\rangle =\left\langle \nabla _{X}(\nabla
f),{Y}\right\rangle ,\;\;X,Y\in {\mathfrak{X}}(M).
\end{equation*}

Note that $\mathbb{L}_{k}$ is the linearized operator of the first variation
of the $(k+1)$-th mean curvature arising from the normal variations of $M$
for $k=2,3,...,n-1$. In particular, $\mathbb{L}_{0}=-\Delta $
is called the Laplace--Beltrami operator (i.e., Laplacian), $\mathbb{L}%
_{1}=\square $ is called the Cheng--Yau operator introduced by Cheng and Yau 
\cite{Cheng-Yau}.

The second order differential operator $\mathbb{L}_{k}$
associated with the Newton transformation $\mathcal{P}_{k}$ is given by 
$\mathbb{L}_{k}\left( f\right) =\mathrm{div}\left( \mathcal{P}%
_{k}\left( \nabla ^{2}f\right) \right) $, \cite{AliasGurbuzGeomDed2006}.
From this equation, the following appears
\begin{equation}
\mathbb{L}_{k}=\sum\limits_{i}r_{k}^{i}\left( e_{i}e_{i}-\nabla
_{e_{i}}e_{i}\right),  \label{LkAltDef}
\end{equation}%
and the equation 
\begin{equation}
\mathbb{L}_{k}x=s_{k}\mathcal{G}  \label{Lkx}
\end{equation}%
is obtained in \cite{AliasGurbuzGeomDed2006}, where $\mathcal{G}$ is the
Gauss map of $M$ which assigns every point of $M$ into the unit normal
vector associated with the oriantation of $M$. Furthermore, bu using %
\eqref{LkAltDef}, one can obtain 
\begin{equation}
\mathbb{L}_{k}\mathcal{G}=-\nabla
s_{k+1}-(s_{1}s_{k+1}-(k+2)s_{k+2})\mathcal{G}.  \label{LkG}
\end{equation}
To thoroughly investigate the $\mathbb{L}_{k}$ operators, we suggest readers
to refer to the comprehensive research conducted by Chen et al. in their
work \cite{CGYH}.

\section{Rotational hypersurfaces in Euclidean $n$-spaces}

We note that the definition of rotational hypersurfaces in Riemannian space
forms were defined in \cite{CDRotHyp}. A rotational hypersurface $M\subset 
\mathbb{E}^{n}$ generated by a curve $\mathcal{C}$ around an axis $\ell$ that does not meet $\mathcal{C}$ is obtained by taking the orbit of $\mathcal{C}$ under those orthogonal transformations of $\mathbb{E}^{n}$ that
leave $\ell$ pointwise fixed. See \cite[Remark 2.3]{CDRotHyp} for
details.

Let $\mathcal{C}$ be the curve of $\mathbb{E}^{n}$ parametrized by $\gamma
=\gamma (r)$: 
\begin{equation}
\gamma (r)=\left( f(r),0,\hdots,0,\varphi \left( r\right) \right) .
\label{0a}
\end{equation}%
An orthogonal
transformations of $\mathbb{E}^{n}$ that leaves $\ell
$ pointwise fixed
has the following $n\times n$ matrix form $\mathcal{Z}=\mathcal{Z}(\theta
_{1},\theta _{2},\hdots,\theta _{n-2})$: 
\begin{equation}
\mathcal{Z}=\left( 
\begin{array}{cccccccc}
\prod\limits_{i=1}^{n-2}C_{i} & -S_{1} & -C_{1}S_{2} & -\prod%
\limits_{i=1}^{2}C_{i}S_{3} & \cdots & -\prod\limits_{i=1}^{n-4}C_{i}S_{n-3}
& -\prod\limits_{i=1}^{n-3}C_{i}S_{n-2} & 0 \\ 
S_{1}\prod\limits_{i=2}^{n-2}C_{i} & C_{1} & -S_{1}S_{2} & 
-S_{1}\prod\limits_{i=2}^{2}C_{i}S_{3} & \cdots & -S_{1}\prod%
\limits_{i=2}^{n-4}C_{i}S_{n-3} & -S_{1}\prod\limits_{i=2}^{n-3}C_{i}S_{n-2}
& 0 \\ 
S_{2}\prod\limits_{i=3}^{n-2}C_{i} & 0 & C_{2} & -S_{2}S_{3} & \cdots & 
-S_{2}\prod\limits_{i=3}^{n-4}C_{i}S_{n-3} & -S_{2}\prod%
\limits_{i=3}^{n-3}C_{i}S_{n-2} & 0 \\ 
S_{3}\prod\limits_{i=4}^{n-2}C_{i} & \vdots & 0 & C_{3} & \cdots & \vdots & 
-S_{3}\prod\limits_{i=4}^{n-3}C_{i}S_{n-2} & 0 \\ 
\vdots & 0 & \vdots & 0 & \ddots & -S_{n-4}S_{n-3} & \vdots & 0 \\ 
S_{n-3}C_{n-2} & 0 & 0 & \vdots & \cdots & C_{n-3} & -S_{n-3}S_{n-2} & \vdots
\\ 
S_{n-2} & 0 & 0 & 0 & \cdots & 0 & C_{n-2} & 0 \\ 
0 & 0 & 0 & 0 & \cdots & 0 & 0 & 1%
\end{array}%
\right) ,  \label{Z}
\end{equation}%
where $\mathcal{Z}\in SO\left( n\right) ,$ $\mathcal{Z}\ell =\ell ,$ $\det 
\mathcal{Z}=1,$ $C_{i}=\cos \theta _{i},$ $S_{i}=\sin \theta _{i},$ $\theta
_{i}\in {\mathbb{R}},$ $\ell $ denotes the rotation axis $(0,0,\hdots%
,0,1)^{T}$.

Next, we define the parametrization of the rotational hypersurfaces family,
generated by a curve $\gamma $, represented by Eq. $\left( \ref{0a}\right) $%
, around axis $\ell $. This parametrization is given by 
\begin{equation}
\mathfrak{x}=\mathcal{Z}\text{\textperiodcentered }\gamma ^{T}.  \label{1}
\end{equation}

In the subsequent sections of this paper, a vector or matrix $(v_{1},v_{2},%
\hdots,v_{n})$ is considered equivalent to its transpose, denoted as the
column vector or column matrix $(v_{1},v_{2},\hdots,v_{n})^{T}$.

Let $\mathfrak{x}=\mathfrak{x}(r,\theta _{1},\hdots,\theta _{n-2})$ be a
parametric representation and isometric immersion of a hypersurfaces family $%
M$ in the Euclidean space $\mathbb{E}^{n}$. 

The vector product of $\mathbf{v}%
_{1}=(v_{1}^{1},v_{1}^{2},\hdots,v_{1}^{n})$,$\hdots$, $\mathbf{v}%
_{n-1}=(v_{n-1}^{1},v_{n-1}^{2},\hdots,v_{n-1}^{n})$ in $\mathbb{E}^{n}$ is
defined by%
\begin{equation*}
\mathbf{v}_{1}\times \mathbf{v}_{2}\times \hdots\times \mathbf{v}_{n-1}=\det
\left( 
\begin{array}{cccc}
e_{1} & e_{2} & \cdots & e_{n} \\ 
v_{1}^{1} & v_{1}^{2} & \cdots & v_{1}^{n} \\ 
v_{2}^{1} & v_{2}^{2} & \cdots & v_{2}^{n} \\ 
\vdots & \vdots & \ddots & \vdots \\ 
v_{n-1}^{1} & v_{n-1}^{2} & \cdots & v_{n-1}^{n}%
\end{array}%
\right) _{n\times n}.
\end{equation*}%
The first and the second fundamental form matrices of hypersurface $%
\mathfrak{x}$ in $\mathbb{E}^{n}$ are described, respectively, by $\mathbf{I}%
=\left( 
\begin{array}{c}
g_{ij}%
\end{array}%
\right) _{\left( n-1\right) \times \left( n-1\right) }$ and $\mathbf{II}%
=\left( 
\begin{array}{c}
h_{ij}%
\end{array}%
\right) _{\left( n-1\right) \times \left( n-1\right) }.$ Here, $1\leq
i,j\leq n-1$, $g_{11}=\left\langle \mathfrak{x}_{r},\mathfrak{x}%
_{r}\right\rangle ,$ $g_{12}=\left\langle \mathfrak{x}_{r},\mathfrak{x}%
_{\theta _{1}}\right\rangle ,\hdots,$ $g_{\left( n-1\right) \left(
n-1\right) }=\left\langle \mathfrak{x}_{\theta _{n-2}},\mathfrak{x}_{\theta
_{n-2}}\right\rangle ,$ $h_{11}=\left\langle \mathfrak{x}_{rr},\mathcal{G}%
\right\rangle ,$ $h_{12}=\left\langle \mathfrak{x}_{r\theta _{1}},\mathcal{G}%
\right\rangle ,\hdots,$ $h_{\left( n-1\right) \left( n-1\right)
}=\left\langle \mathfrak{x}_{\theta _{n-2}\theta _{n-2}},\mathcal{G}%
\right\rangle $, the partial differentials are $\mathfrak{x}_{r}=\frac{%
\partial \mathfrak{x}}{\partial r}$, $\mathfrak{x}_{r\theta _{n-2}}=\frac{%
\partial ^{2}\mathfrak{x}}{\partial r\partial \theta _{n-2}}$, etc., and the
Gauss map (i.e., the unit normal vector) of hypersurface $\mathfrak{x}$ is
expressed by 
\begin{equation}
\mathcal{G}=\frac{\mathfrak{x}_{r}\times \mathfrak{x}_{\theta _{1}}\times %
\hdots\times \mathfrak{x}_{\theta _{n-2}}}{\left\Vert \mathfrak{x}_{r}\times 
\mathfrak{x}_{\theta _{1}}\times \hdots\times \mathfrak{x}_{\theta
_{n-2}}\right\Vert }.  \label{2}
\end{equation}%
$\left( g_{ij}\right) ^{-1}\left( h_{ij}\right) $ gives the shape operator
matrix $\mathbf{A}.$ Then, the formulas of the mean curvature and the
Gauss--Kronecker curvature are, respectively, determined by%
\begin{equation}
H=\frac{1}{n-1}\text{tr}\left( \mathbf{A}\right) ,\text{ and }K=\det (%
\mathbf{A})=\frac{\det \mathbf{II}}{\det \mathbf{I}}.  \label{4}
\end{equation}

\section{Curvatures of the rotational hypersurfaces family in $\mathbb{E}%
^{n} $}

By considering Eq. $\left( \ref{1}\right) $, the following  parametric form
of the rotational hypersurfaces family $\mathfrak{x}=\mathfrak{x}(r,\theta
_{1},\theta _{2},\hdots,\theta _{n-2})$ appears
\begin{eqnarray}
\mathfrak{x} &=&\left( f(r)\prod\limits_{i=1}^{n-2}\cos \theta
_{i},f(r)\sin \theta _{1}\prod\limits_{i=2}^{n-2}\cos \theta _{i},\right.
\label{6} \\
&&\text{ \ \ \ \ \ \ \ \ \ \ \ }\left. f(r)\sin \theta
_{2}\prod\limits_{i=3}^{n-2}\cos \theta _{i},\hdots,f(r)\sin \theta
_{n-3}\cos \theta _{n-2},f(r)\sin \theta _{n-2},\varphi (r)\right)  \notag \\
&=&\left( x_{1},x_{2},x_{3},\hdots,x_{n-2},x_{n-1},x_{n}\right) ,  \notag
\end{eqnarray}%
where $f=f(r),$ $r,\theta _{i}\in \mathbb{R\setminus \{}0\mathbb{\}}.$
By taking the first differentials with respect to $r,\theta _{1},\theta _{2},%
\hdots,\theta _{n-2}$ of the rotational hypersurfaces family given by Eq. $\left(\ref{6}\right)$, the following first fundamental form matrix occurs
\begin{equation*}
\mathbf{I}=\text{diag}\left( 
\begin{array}{cccccc}
f^{\prime 2}+\varphi ^{\prime 2}, & f^{2}\prod\limits_{i=2}^{n-2}\cos
^{2}\theta _{i}, & f^{2}\prod\limits_{i=3}^{n-2}\cos ^{2}\theta _{i}, & %
\hdots, & f^{2}\cos ^{2}\theta _{n-2}, & f^{2}%
\end{array}%
\right).
\end{equation*}%
Then,
\begin{equation*}
\det \mathbf{I}=\left( f^{2}\right) ^{n-2}(f^{\prime 2}+\varphi ^{\prime
2})\left( \prod\limits_{i=2}^{n-2}\cos \theta _{i}\right) \left(
\prod\limits_{i=3}^{n-2}\cos ^{2}\theta _{i}\right) \hdots\left(
\prod\limits_{i=n-2}^{n-2}\cos ^{2}\theta _{i}\right) ,
\end{equation*}%
where $f=f(r),$ $f^{\prime }=\frac{df}{dr},$ $\varphi =\varphi (r),$ $%
\varphi ^{\prime }=\frac{d\varphi }{dr}.$ By taking the second differentials of $%
f$ and $\theta _{i}$ of the rotational hypersurfaces family determined by Eq. $%
\left(\ref{6}\right) $, the following second fundamental form matrix appears
\begin{equation*}
\mathbf{II}=\varepsilon \text{ diag}\left( 
\begin{array}{cccccc}
\frac{f^{\prime }\varphi ^{\prime \prime }-f^{\prime \prime }\varphi
^{\prime }}{(f^{\prime 2}+\varphi ^{\prime 2})^{1/2}}, & \frac{f\varphi
^{\prime }\prod\limits_{i=2}^{n-2}\cos ^{2}\theta _{i}}{(f^{\prime
2}+\varphi ^{\prime 2})^{1/2}}, & \frac{f\varphi ^{\prime
}\prod\limits_{i=3}^{n-2}\cos ^{2}\theta _{i}}{(f^{\prime 2}+\varphi
^{\prime 2})^{1/2}}, & \hdots, & \frac{f\varphi ^{\prime }\cos ^{2}\theta
_{n-2}}{(f^{\prime 2}+\varphi ^{\prime 2})^{1/2}}, & \frac{f\varphi ^{\prime
}}{(f^{\prime 2}+\varphi ^{\prime 2})^{1/2}}%
\end{array}%
\right).
\end{equation*}%
Then,
\begin{equation*}
\det \mathbf{II}=\varepsilon \frac{f^{n-2}\left( f^{\prime }\varphi ^{\prime
\prime }-f^{\prime \prime }\varphi \right) \left( \varphi ^{\prime }\right)
^{n-2}\left( \prod\limits_{i=2}^{n-2}\cos ^{2}\theta _{i}\right) \left(
\prod\limits_{i=3}^{n-2}\cos ^{2}\theta _{i}\right) \hdots\left(
\prod\limits_{i=n-2}^{n-2}\cos ^{2}\theta _{i}\right) }{(f^{\prime
2}+\varphi ^{\prime 2})^{\left( n-1\right) /2}}.
\end{equation*}%
Here, 
\begin{equation}
\varepsilon =\left\{ 
\begin{array}{cll}
-1 & \text{if} & n\text{ odd integer,} \\ 
1 & \text{if} & n\text{ even integer,}%
\end{array}%
\right.  \label{eq.epsilon}
\end{equation}%
for integers $n\geq 3.$ The Gauss map of the rotational hypersurfaces family
described by Eq. $\left(\ref{6}\right)$ is determined by%
\begin{equation*}
\mathcal{G}=\frac{\varepsilon }{(f^{\prime 2}+\varphi ^{\prime 2})^{1/2}}%
\left( \varphi ^{\prime }\prod\limits_{i=1}^{n-2}C_{i},\varphi ^{\prime
}S_{1}\prod\limits_{i=2}^{n-2}C_{i},\varphi ^{\prime
}S_{2}\prod\limits_{i=3}^{n-2}C_{i},\hdots,\varphi ^{\prime
}S_{n-3}C_{n-2},\varphi ^{\prime }S_{n-2},-f^{\prime }\right) .
\end{equation*}%
The shape operator of the family determined by Eq. $\left( \ref{6}\right) $
is given by 
\begin{equation*}
\mathbf{A}=\text{diag}\left( 
\begin{array}{cccc}
\kappa _{1}, & \kappa _{2}, & \hdots, & \kappa _{n-1}%
\end{array}%
\right) .
\end{equation*}%
Here, the principal curvatures of the rotational family defined by Eq. $\left(\ref{6}\right)$ are described by $\kappa _{1}=\varepsilon \dfrac{%
f^{\prime }\varphi ^{\prime \prime }-f^{\prime \prime }\varphi ^{\prime }}{%
(f^{\prime 2}+\varphi ^{\prime 2})^{3/2}},$ $\kappa _{2}=\varepsilon \dfrac{%
\varphi ^{\prime }}{f(f^{\prime 2}+\varphi ^{\prime 2})^{1/2}}=\kappa _{3}=%
\hdots=\kappa _{n-1}.$ Then, the mean curvature is given by 
\begin{equation}
H=\varepsilon \frac{ff^{\prime }\varphi ^{\prime \prime }+\left( n-2\right)
\varphi ^{\prime 3}+\left( \left( n-2\right) f^{\prime 2}-ff^{\prime \prime
}\right) \varphi ^{\prime }}{\left( n-1\right) f\left( f^{\prime 2}+\varphi
^{\prime 2}\right) ^{3/2}},  \label{meanH}
\end{equation}%
and the Gauss--Kronecker curvature is determined by%
\begin{equation}
K=\varepsilon \frac{\left( f^{\prime }\varphi ^{\prime \prime }-f^{\prime
\prime }\varphi ^{\prime }\right) \left( \varphi ^{\prime }\right) ^{n-2}}{%
f^{n-2}(f^{\prime 2}+\varphi ^{\prime 2})^{n/2}}.  \label{GaussK}
\end{equation}

Therefore, the following theorems come into view.

\begin{theorem} 
\textit{The rotational hypersurfaces family parametrized
by }
\begin{eqnarray*}
\mathfrak{x} &=&\left( f(r)\prod\limits_{i=1}^{n-2}\cos \theta
_{i},f(r)\sin \theta _{1}\prod\limits_{i=2}^{n-2}\cos \theta _{i},\right. \\
&&\text{ \ \ \ \ \ \ \ \ \ \ \ }\left. f(r)\sin \theta
_{2}\prod\limits_{i=3}^{n-2}\cos \theta _{i},\hdots,f(r)\sin \theta
_{n-3}\cos \theta _{n-2},f(r)\sin \theta _{n-2},\varphi (r)\right)
\end{eqnarray*}%
\textit{is minimal }$($i.e., $H=0)$ \textit{if and only if the following
holds }%
\begin{equation*}
\varphi (r)=\pm \int \dfrac{f^{2}(r)f^{\prime }(r)}{\left(
	c_{1}f^{2n}(r)-f^{4}(r)\right) ^{1/2}}dr+c_{2},
	\end{equation*}
	\textit{where} $c_{1},c_{2}\in \mathbb{R}$\textit{,} \textit{and }$%
	(n-1)f(f^{\prime 2}+\varphi ^{\prime 2})^{3/2}\neq 0.$
\end{theorem} 

\begin{proof}
When considering the rotational hypersurface determined by Eq. $%
\left( \ref{6}\right) $ under the condition of minimality, i.e., when $H=0$ 
\textit{described by Eq.} $\left(\ref{meanH}\right)$, we arrive at a
second-order differential equation $ff^{\prime }\varphi ^{\prime \prime
}+\left( n-2\right) \varphi ^{\prime 3}+\left( \left( n-2\right) f^{\prime
2}-ff^{\prime \prime }\right) \varphi ^{\prime }=0.$ The solutions for $\varphi (r)$ could not be obtained using classical methods. 

However, $\varphi (r)$ is successfully obtained as indicated in the theorem, using Maple program codes:
\begin{eqnarray*}
	&>&dsolve(f(r)\ast diff(f(r),r)\ast diff(\varphi (r),r,r)+(n-2)\ast
	diff(\varphi (r),r)\symbol{94}3 \\
	&&+((n-2)\ast diff(f(r),r)\symbol{94}2-f(r)\ast diff(f(x),r,r))\ast
	diff(\varphi (r),r) \\
	&=&0,\varphi (r));
\end{eqnarray*}
\end{proof}

\begin{theorem} 
\textit{\ The rotational hypersurfaces family defined by
Eq.} $\left( \ref{6}\right) $\textit{\ has zero Gauss--Kronecker curvature }$%
($i.e., $K=0)$\textit{\ if and only if }\begin{equation*}
	\varphi (r)=c_{1}\ \text{\ \ or \ }\varphi (r)=c_{1}f(r)+c_{2},
\end{equation*}
\textit{where }$f^{n-2}(f^{\prime 2}+\varphi ^{\prime 2})^{n/2}\neq 0,$ $%
c_{1},c_{2}\in \mathbb{R}.$
\end{theorem} 

\begin{proof}
When $K=0$ \textit{described by Eq. }$\left( \ref{GaussK}\right) ,$
the second-order differential equation $f^{\prime }\left( \varphi ^{\prime
}\right) ^{n-2}\varphi ^{\prime \prime }-f^{\prime \prime }\left( \varphi
^{\prime }\right) ^{n-1}=0$ transforms to $\dfrac{\varphi ^{\prime \prime }}{%
\varphi ^{\prime }}=\dfrac{f^{\prime \prime }}{f^{\prime }}.$ Then, we
obtain solutions $\varphi (r)$. 

We also find solution $\varphi (r)=f(r)/c_{1}+c_{2}$ using a Maple software
program codes: 
\begin{eqnarray*}
	&>&PDEtools[declare](f(r),\varphi (r),prime=r); \\
	&>&dsolve(diff(f(r),r)\ast diff(\varphi (r),r)\symbol{94}(n-2)\ast
	diff(\varphi (r),r,r) \\
	&&-diff(f(r),r,r)\ast diff(\varphi (r),r)\symbol{94}(n-1) \\
	&=&0,\varphi (r));
\end{eqnarray*}
\end{proof}

\begin{corollary} 
\textit{Substituting }$\varphi =c=const.$\textit{\
into the Eqs. described by }$\left( \ref{meanH}\right) $ \textit{and} $%
\left( \ref{GaussK}\right) ,$ \textit{we find the curvatures }$H=0,$\ $K=0.$ 
\textit{Then, the hypersurface is minimal.}
\end{corollary} 

\begin{corollary} 
\textit{Substituting }$f=c=const.$ \textit{into the
Eqs. determined by }$\left( \ref{meanH}\right) $ \textit{and} $\left( \ref%
{GaussK}\right) ,$ \textit{we obtain the curvatures }$H=-\varepsilon \dfrac{%
n-2}{\left( n-1\right) c},$ $K=0,$ \textit{where} $c\neq 0.$ \textit{Hence,
the hypersurface has CMC.}
\end{corollary}

\section{Gauss map and $\mathbb{L}_{n-3}$ operator in $\mathbb{E}^{n}$}

We consider the parametric curve $\gamma =\gamma \left( r\right) $ defined
by Eq. $\left( \ref{0a}\right) $ as a curve with unit speed. The elements of
the adapted frame field $\left\{ e_{1},e_{2},\hdots,e_{n-1},\mathcal{G}%
\right\} $ of the rotational hypersurfaces family $\mathfrak{x}$ determined
by Eq\textit{.} $\left( \ref{6}\right) $ are given by%
\begin{eqnarray}
e_{1} &=&\mathfrak{x}_{r}  \label{e_1} \\
&=&\left( f^{\prime }\prod\limits_{i=1}^{n-2}\cos \theta _{i},f^{\prime
}\sin \theta _{1}\prod\limits_{i=2}^{n-2}\cos \theta _{i},\right.   \notag
\\
&&\text{ \ \ \ \ \ \ \ \ \ \ \ }\left. f^{\prime }\sin \theta
_{2}\prod\limits_{i=3}^{n-2}\cos \theta _{i},\hdots,f^{\prime }\sin \theta
_{n-3}\cos \theta _{n-2},f^{\prime }\sin \theta _{n-2},\varphi ^{\prime
}\right) ,  \notag
\end{eqnarray}%
\begin{eqnarray}
e_{2} &=&\frac{\mathfrak{x}_{\theta _{1}}}{f}  \label{e_2} \\
&=&\left( 
\begin{array}{c}
-\sin \theta _{1}\prod\limits_{i=2}^{n-2}\cos \theta _{i},\cos \theta
_{1}\prod\limits_{i=2}^{n-2}\cos \theta _{i},0,\hdots,0%
\end{array}%
\right) ,  \notag
\end{eqnarray}%
\begin{equation*}
\vdots 
\end{equation*}%
\begin{eqnarray}
e_{n-1} &=&\frac{\mathfrak{x}_{\theta _{n-2}}}{f}  \label{e_n-1} \\
&=&\left( -\sin \theta _{n-2}\prod\limits_{i=1}^{n-3}\cos \theta _{i},-\sin
\theta _{n-2}\prod\limits_{i=2}^{n-3}\cos \theta _{i},\right.   \notag \\
&&\text{ \ \ \ \ \ \ \ \ \ \ \ }\left. -\sin \theta
_{n-2}\prod\limits_{i=3}^{n-3}\cos \theta _{i},\hdots,-\sin \theta
_{n-3}\sin \theta _{n-2},\cos \theta _{n-2},0\right) ,  \notag
\end{eqnarray}%
\begin{eqnarray}
\mathcal{G} &=&e_{1}\times e_{2}\times \hdots\times e_{n-1}  \label{G} \\
&=&\varepsilon \left( \varphi ^{\prime }\prod\limits_{i=1}^{n-2}\cos \theta
_{i},\varphi ^{\prime }\sin \theta _{1}\prod\limits_{i=2}^{n-2}\cos \theta
_{i},\right.   \notag \\
&&\text{ \ \ \ \ \ \ \ \ \ \ \ }\left. \varphi ^{\prime }\sin \theta
_{2}\prod\limits_{i=3}^{n-2}\cos \theta _{i},\hdots,\varphi ^{\prime }\sin
\theta _{n-3}\cos \theta _{n-2},\varphi ^{\prime }\sin \theta
_{n-2},-f^{\prime }\right) .  \notag
\end{eqnarray}%
$\bigskip $
The principal curvatures $\kappa _{i}$ $\left( 1\leq i\leq n-1\right) $ of $%
\mathfrak{x}$ with respect to the Gauss map $\mathcal{G}$ are determined by
\begin{eqnarray*}
\kappa _{1} &=&\left\langle \mathbf{A}\left( e_{1}\right)
,e_{1}\right\rangle =-\varepsilon \left( f^{\prime \prime }\varphi ^{\prime
}-f^{\prime }\varphi ^{\prime \prime }\right) , \\
\kappa _{j} &=&\left\langle \mathbf{A}\left( e_{j}\right)
,e_{j}\right\rangle =-\varepsilon \frac{\varphi ^{\prime }}{f},\text{ \ \ }%
\left( 2\leq j\leq n-1\right) .
\end{eqnarray*}
Since $\gamma $ denotes a unit speed curve, there exists a smooth function $%
R=R(r)$ such that $f^{\prime }=\cos R$ and $\varphi ^{\prime }=\sin R$.
Then, the following is presented.

\begin{lemma} 
\textit{Consider}\ \textit{an\ oriented\ hypersurface\ }$%
\mathfrak{x}$\textit{\ in\ the\ Euclidean\ space}\ $\mathbb{E}^{n}$,\ 
\textit{with\ its\ mean\ curvature}\ $H$\textit{\ and\ Gauss--Kronecker\
curvature\ }$K$.\ \textit{Then,\ the\ Gauss\ map\ of}\ $\mathfrak{x}$\textit{%
,\ for\ integers}\ $n\geq 3$\textit{,}\ \textit{satisfies\ the\ relation}%
\begin{equation}
\mathbb{L}_{n-3}\mathcal{G}=-\nabla s_{n-2}-(s_{1}s_{n-2}-\left( n-1\right)
s_{n-1})\mathcal{G},  \label{L_(n-3) G}
\end{equation}%
\textit{where}%
\begin{equation}
\begin{array}{l}
s_{1}=\varepsilon \left( \dbinom{n-2}{0}R^{\prime }-\dbinom{n-2}{1}\wp
\right) , \\ 
s_{n-2}=\varepsilon \left( \dbinom{n-2}{n-3}R^{\prime }\wp ^{n-3}-\dbinom{n-2%
}{n-2}\wp ^{n-2}\right) , \\ 
s_{n-1}=\varepsilon \dbinom{n-2}{n-2}R^{\prime }\wp ^{n-2},%
\end{array}
\label{s1s2s3...s(n-1)}
\end{equation}%
$s_{1}=\left( n-1\right) H,$ $s_{n-1}=K,$ $\wp =\frac{\sin R}{f},$ $\gamma
(r)=\left( f(r),0,\hdots,0,\varphi \left( r\right) \right) $\ \textit{%
represents\ a\ unit\ speed\ curve},\ \textit{and}\ $R=R(r)$\ \textit{%
denotes\ a\ smooth\ function,} $f^{\prime }=\cos R$, $\varphi
^{\prime }=\sin R$.
\end{lemma} 

Before proving Lemma 1, let's provide the following examples.

\bigskip

\textit{Example 1.} \textit{Flat hypersurfaces} \textit{satisfy} $\mathbb{L}%
_{n-3}\mathcal{G}=\mathcal{AG}$\textit{,} \textit{where} $\mathcal{A}$ 
\textit{denotes a} $n\times n$\textit{\ matrix, derived from} $\mathbb{L}%
_{n-3}\mathcal{G}=-\nabla s_{n-2}-s_{1}s_{n-2}\mathcal{G}.$

\bigskip

\textit{Example 2.} \textit{Hyperspheres}\textbf{\ }$\sum\limits_{i=1}^{n}\left(
x_{i}-p_{i}\right) ^{2}=r^{2}\ $\textit{have} $\mathcal{G}=\frac{1}{r}\left(
x_{1}-p_{1},\hdots,x_{n}-p_{n}\right) $. \textit{Then,} \textit{hyperspheres
hold }$\mathbb{L}_{n-3}\mathcal{G}=\mathcal{AG}$\textit{,} \textit{with }$%
\mathcal{A}=\varepsilon \frac{1}{r^{n}}\mathcal{I}_{n}$\textit{,} \textit{%
where} $\mathcal{I}_{n}$ \textit{represents an} \textit{identity matrix.}

\bigskip

\textit{Example 3.} \textit{Minimal hypersurfaces satisfy} $\mathbb{L}%
_{n-3}\mathcal{G}=-\nabla s_{n-2}+\left( n-1\right) s_{n-1}\mathcal{G}$.

\bigskip

Now, let's  prove  Lemma 1.

\begin{proof}
Taking $f^{\prime }=\cos R$ and $\varphi
^{\prime }=\sin R$, the gradient $\nabla s_{n-2}$ of $\mathfrak{x}$\ is given by
\begin{equation}
\nabla s_{n-2}=s_{n-2}^{\prime }e_{1},  \label{gradient s _(n-2)}
\end{equation}%
where%
\begin{eqnarray}
s_{n-2}^{\prime }(r) &=&\varepsilon \frac{1}{f^{n-1}}\left( n-2\right)
\left( f^{2}R^{\prime \prime }\sin R+\left( n-3\right) f^{2}R^{\prime 2}\cos
R\right.  \label{s' _(n-2)} \\
&&\left. -\left( n-4\right) fR^{\prime }\sin R\cos R-\sin ^{2}R\cos R\right)
\sin ^{n-4}R.  \notag
\end{eqnarray}%
Next, we suppose that the Gauss map $\mathcal{G}$  given by Eq.\textit{\ }$\left( \ref{G}%
\right)$ of the rotational
hypersurface $\mathfrak{x}$ satisfies the following
\begin{equation}
\mathbb{L}_{n-3}\mathcal{G}=\mathcal{AG}.  \label{L_(n-3) G=AG}
\end{equation}%
 with a $n\times n$ matrix $\mathcal{A}%
=\left( a_{ij}\right) $. Substituting $f^{\prime }=\cos R$ and $\varphi ^{\prime }=\sin R$ again,
into the Gauss map $\mathcal{G}$ determined by Eq.\textit{\ }$\left( \ref{G}%
\right) $, the following holds
\begin{eqnarray}
\mathcal{G} &=&\varepsilon \left( \sin R\prod\limits_{i=1}^{n-2}\cos \theta
_{i},\sin R\sin \theta _{1}\prod\limits_{i=2}^{n-2}\cos \theta _{i},\right.
\label{newG} \\
&&\left. \sin R\sin \theta _{2}\prod\limits_{i=3}^{n-2}\cos \theta _{i},%
\hdots,\sin R\sin \theta _{n-3}\cos \theta _{n-2},\sin R\sin \theta
_{n-2},-\cos R\right) .  \notag
\end{eqnarray}%
Consequently, deducing from Eqs. $\left( \ref{s1s2s3...s(n-1)}\right) ,$ $%
\left( \ref{L_(n-3) G=AG}\right) $ and $\left( \ref{newG}\right) ,$ we infer
the following system 
\begin{equation}
\begin{array}{l}
\left( -s_{n-2}^{\prime }\cos R-(s_{1}s_{n-2}-\left( n-1\right) s_{n-1})\sin
R\right) \cos \theta _{1}\hdots\cos \theta _{n-2} \\ 
=\varepsilon (a_{11}\sin R\cos \theta _{1}\hdots\cos \theta
_{n-2}+a_{12}\sin R\sin \theta _{1}\cos \theta _{2}\hdots\cos \theta _{n-2}
\\ 
\text{ \ \ }+\hdots+a_{1\left( n-2\right) }\sin R\sin \theta _{n-3}\cos
\theta _{n-2}+a_{1\left( n-1\right) }\sin R\sin \theta _{n-2}-a_{1n}\cos R),
\\ 
\end{array}
\label{L_(n-3) G=AG 1}
\end{equation}%
\begin{equation}
\begin{array}{l}
\left( -s_{n-2}^{\prime }\cos R-(s_{1}s_{n-2}-\left( n-1\right) s_{n-1})\sin
R\right) \sin \theta _{1}\cos \theta _{2}\hdots\cos \theta _{n-2} \\ 
=\varepsilon (a_{21}\sin R\cos \theta _{1}\hdots\cos \theta
_{n-2}+a_{22}\sin R\sin \theta _{1}\cos \theta _{2}\hdots\cos \theta _{n-2}
\\ 
\text{ \ \ }+\hdots+a_{2\left( n-2\right) }\sin R\sin \theta _{n-3}\cos
\theta _{n-2}+a_{2\left( n-1\right) }\sin R\sin \theta _{n-2}-a_{2n}\cos R),%
\end{array}
\label{L_(n-3) G=AG 2}
\end{equation}%
\begin{equation*}
\vdots
\end{equation*}%
\begin{equation}
\begin{array}{l}
\left( -s_{n-2}^{\prime }\cos R-(s_{1}s_{n-2}-\left( n-1\right) s_{n-1})\sin
R\right) \sin \theta _{n-2} \\ 
=\varepsilon (a_{\left( n-1\right) 1}\sin R\cos \theta _{1}\hdots\cos \theta
_{n-2}+a_{\left( n-1\right) 2}\sin R\sin \theta _{1}\cos \theta _{2}\hdots%
\cos \theta _{n-2} \\ 
\text{ \ \ }+\hdots+a_{\left( n-1\right) \left( n-2\right) }\sin R\sin
\theta _{n-3}\cos \theta _{n-2}+a_{\left( n-1\right) \left( n-1\right) }\sin
R\sin \theta _{n-2}-a_{\left( n-1\right) n}\cos R),%
\end{array}
\label{L_(n-3) G=AG  n-1}
\end{equation}%
\begin{equation}
\begin{array}{l}
-s_{n-2}^{\prime }\sin R-(s_{1}s_{n-2}-\left( n-1\right) s_{n-1})\cos R \\ 
=\varepsilon (a_{n1}\sin R\cos \theta _{1}\hdots\cos \theta
_{n-2}+a_{n2}\sin R\sin \theta _{1}\cos \theta _{2}\hdots\cos \theta _{n-2}
\\ 
\text{ \ \ }+\hdots+a_{n\left( n-2\right) }\sin R\sin \theta _{n-3}\cos
\theta _{n-2}+a_{n\left( n-1\right) }\sin R\sin \theta _{n-2}-a_{nn}\cos R).%
\end{array}
\label{L_(n-3) G=AG  n}
\end{equation}%
Next, we suppose that $J=\left\{ r\in I\mid R^{\prime }\left( r\right) \neq
0\right\} $ is nonempty set. Then, $R(I)$ contains an interval, and obtain
from $\left( \ref{L_(n-3) G=AG 1}\right) $-$\left( \ref{L_(n-3) G=AG n}%
\right) $ that $a_{11}=a_{22}=\hdots=a_{\left( n-1\right) \left( n-1\right)
} $ and $a_{ij}=0$ when $i\neq j.$ Hence, the following hold $\mathcal{A}=$ diag$\left( 
\underset{n}{\underbrace{\eta ,\hdots,\eta ,\phi }}\right) ,$%
\begin{equation}
-s_{n-2}^{\prime }\cos R-(s_{1}s_{n-2}-\left( n-1\right) s_{n-1})\sin
R=\varepsilon \eta \sin R,  \label{Eq1}
\end{equation}%
and%
\begin{equation}
-s_{n-2}^{\prime }\sin R-(s_{1}s_{n-2}-\left( n-1\right) s_{n-1})\cos
R=-\varepsilon \phi \cos R.  \label{Eq2}
\end{equation}%
Eqs. $\left( \ref{Eq1}\right) $ and $\left( \ref{Eq2}\right) $ are equivalent to
the following%
\begin{equation}
s_{n-2}^{\prime }=\lambda \sin R\cos R,  \label{Eq3}
\end{equation}%
\begin{equation}
s_{1}s_{n-2}-\left( n-1\right) s_{n-1}=\varepsilon \left( \lambda \sin
^{2}R+\phi \right) ,  \label{Eq4}
\end{equation}%
where $\lambda =\phi -\eta .$
\end{proof}

Therefore, the following is provided.

\begin{lemma} 
\textit{Consider a rotational hypersurface }$\mathfrak{x}$%
\textit{\ given by} $\left( \ref{6}\right) $ \textit{in Euclidean space }$%
\mathbb{E}^{n}$\textit{,} \textit{with the set} $J=\left\{ r\in I\mid
R^{\prime }\left( r\right) \neq 0\right\} $\textit{. Assume that the Gauss
map} $\mathcal{G}$ \textit{of} $\mathfrak{x}$ \textit{satisfies} $\mathbb{L}%
_{n-3}\mathcal{G}=\mathcal{AG}$\textit{, where} $\mathcal{A}$ \textit{denotes%
} $n\times n$ \textit{matrix with integers} $n\geq 3$\textit{. Then,} $%
\mathcal{A}$ \textit{takes the form} $\eta \mathcal{I}_{n}$\textit{, where }$%
\mathcal{I}_{n}$ \textit{represents the identity matrix.}
\end{lemma} 

\begin{proof}
 We understand from the above findings that $\mathcal{A}$ is a $%
n\times n$ diagonal matrix satisfying $\mathcal{A}=$diag$\left( \eta ,\hdots%
,\eta ,\phi \right) $ for constants $\eta $ and $\phi .$ Then, it follows
from $\left(\ref{s1s2s3...s(n-1)}\right)$, $\left(\ref{s' _(n-2)}\right)$, $\left(\ref{Eq3}\right)$, and $\left(\ref{Eq4}\right)$, the following appear
\begin{eqnarray}
&&\left( n-2\right) \left( f^{2}R^{\prime \prime }\sin R+\left( n-3\right)
f^{2}R^{\prime 2}\cos R\right.   \label{Eq6} \\
&&\left. -\left( n-4\right) fR^{\prime }\sin R\cos R-\sin ^{2}R\cos R\right)
\sin ^{n-4}R  \notag \\
&=&\lambda f^{n-1}\sin R\cos R,  \notag
\end{eqnarray}
and%
\begin{equation}
\left( n-2\right) \left( f^{2}R^{\prime 2}+\sin ^{2}R+\left( n-3\right)
fR^{\prime }\sin R\right) \sin ^{n-3}R=\left( \lambda \sin ^{2}R+\phi
\right) f^{n-1},  \label{Eq7}
\end{equation}%
where $R=R(r).$ Differentiating the both sides of Eq. $\left( \ref{Eq7}\right) $
with respect to $r,$ the following occurs
\begin{equation}
\begin{array}{l}
-2\left( n-2\right) f^{2}R^{\prime }R^{\prime \prime }\sin ^{n-3}R-\left(
n-3\right) \left( n-2\right) fR^{\prime \prime }\sin ^{n-2}R \\ 
-\left( n-3\right) \left( n-2\right) f^{2}R^{\prime 3}\cos R\sin
^{n-4}R-\left( n-2\right) \left( n^{2}-5n+8\right) fR^{\prime 2}\cos R\sin
^{n-3}R \\ 
-2\left( \left( n-2\right) ^{2}\sin ^{n-2}R+\lambda f^{n-1}\sin R\right)
R^{\prime }\cos R-\left( n-1\right) f^{n-2}\left( \lambda \sin ^{2}R+\phi
\right) \cos R=0.
\end{array}   \label{Eq8}
\end{equation}%
Substituting $f^{\prime }=\cos R$ into Eq. $\left(\ref{Eq6}\right)$, the following holds 
\begin{eqnarray}
R^{\prime \prime } &=&\frac{1}{\left( n-2\right) f^{2}\sin ^{n-4}R}[-\left(
n-2\right) \left( n-3\right) f^{2}R^{\prime 2}\cos R\sin ^{n-5}R  \notag \\
&&\left. +\left( n-2\right) \left( n-4\right) fR^{\prime }\cos R\sin
^{n-4}R  +\lambda f^{n-1}\cos R+\left( n-2\right) \cos R\sin ^{n-3}R\right] .
\label{Eq10}
\end{eqnarray}%
Substituting right hand of $R^{\prime \prime }$ into Eq. $\left(\ref{Eq8}%
\right)$, and using $f^{\prime }=\cos R,$ again,
\begin{equation}
\begin{array}{l}
\lbrack \left( n-3\right) \left( n-2\right) f^{3}R^{\prime 3}\sin
^{n-4}R-3\left( n-3\right) \left( n-2\right) f^{2}R^{\prime 2}\sin ^{n-3}R
\\ 
-4\lambda f^{n}R^{\prime }\sin R-\left( n-2\right) \left( n^{2}-5n+10\right)
fR^{\prime }\sin ^{n-2}R \\ 
+\left( -\lambda \left( n-3\right) \sin ^{2}R-\left( n-1\right) \left(
\lambda \sin ^{2}R+\phi \right) \right) f^{n-1}\\ -\left( n-3\right) \left(
n-2\right) \sin ^{n-1}R]\cos R=0.%
\end{array}
\label{Eq11}
\end{equation}%
Taking $R^{\prime 2}$ in $\left( \ref{Eq7}\right)$, and substituting it into 
$\left( \ref{Eq11}\right)$, the following appears
\begin{equation}
R^{\prime }=\dfrac{\alpha _{1}f^{n-1}+\alpha _{2}}{\alpha _{3}f^{n}+\alpha
_{4}f},  \label{Eq13}
\end{equation}%
where%
\begin{equation}
\begin{array}{l}
\alpha _{1}=(\left( n^{2}-7n+14\right) \lambda \sin ^{3}R+\left(
n^{2}-8n+17\right) \phi \sin R)\cos R, \\ 
\alpha _{2}=-\left( n-7\right) \left( n-3\right) \left( n-2\right) \cos
R\sin ^{n}R, \\ 
\alpha _{3}=-(\lambda \left( n+1\right) \sin ^{2}R+\left( n-3\right) \phi
)\cos R, \\ 
\alpha _{4}=\left( n-2\right) (\left( n-3\right) ^{3}-4\left(
n^{2}-6n+10\right) )\cos R\sin ^{n-1}R.%
\end{array}
\label{13a}
\end{equation}%
$\allowbreak $Replace $R^{\prime }$ by Eq. $\left( \ref{Eq7}\right) $
with that by Eq. $\left( \ref{Eq13}\right) .$ Then, 
\begin{equation}
t_{3\left( n-1\right) }f^{3\left( n-1\right) }+t_{2\left( n-1\right)
}f^{2\left( n-1\right) }+t_{n-1}f^{n-1}+t_{0}=0.  \label{Eq14}
\end{equation}%
Here, $t_{i}=t_{i}\left( R\right) $ represent
\begin{equation}
\begin{array}{lll}
t_{3\left( n-1\right) } & = & -\left( \phi +\lambda \sin ^{2}R\right) \alpha
_{3}^{2}, \\ 
t_{2\left( n-1\right) } & = & \left( n-2\right) \alpha _{3}^{2}\sin
^{n-1}R+\left( n-3\right) \left( n-2\right) \alpha _{1}\alpha _{3}\sin
^{n-2}R \\ 
&  & +\left( n-2\right) \alpha _{1}^{2}\sin ^{n-3}R-2\left( \phi +\lambda
\sin ^{2}R\right) \alpha _{3}\alpha _{4}, \\ 
t_{n-1} & = & 2\left( n-2\right) \alpha _{3}\alpha _{4}\sin ^{n-1}R+\left(
n-3\right) \left( n-2\right) \left( \alpha _{1}\alpha _{4}+\alpha _{2}\alpha
_{3}\right) \sin ^{n-2}R \\ 
&  & +2\left( n-2\right) \alpha _{1}\sin ^{n-3}R-\left( \phi +\lambda \sin
^{2}R\right) \alpha _{4}^{2}, \\ 
t_{0} & = & \left( n-2\right) \alpha _{4}^{2}\sin ^{n-1}R+\left( n-3\right)
\left( n-2\right) \alpha _{2}\alpha _{4}\sin ^{n-2}R+\left( n-2\right)
\alpha _{2}^{2}\sin ^{n-3}R.%
\end{array}
\label{14a}
\end{equation}%
Differentiating Eq. $\left( \ref{Eq14}\right) $ with respect to $r,$and
using $f^{\prime }=\cos R,$ $\frac{\partial t_{i}}{\partial r}=\frac{%
\partial t_{i}}{\partial R}\frac{\partial R}{\partial r}=t_{i}^{\prime
}R^{\prime },$ $R^{\prime }(r)$ in Eq. $\left( \ref{Eq13}\right) ,$ the following holds
\begin{equation}
k_{3\left( n-1\right) }f^{3\left( n-1\right) }+k_{2\left( n-1\right)
}f^{2\left( n-1\right) }+k_{n-1}f^{n-1}+k_{0}=0,  \label{Eq15}
\end{equation}%
where $k_{i}=k_{i}\left( R\right) $,%
\begin{equation}
\begin{array}{lll}
k_{3\left( n-1\right) } & = & 3\left( n-1\right) \alpha _{4}t_{3\left(
n-1\right) }+2\left( n-1\right) \alpha _{3}t_{2\left( n-1\right) }+\alpha
_{2}t_{3\left( n-1\right) }^{\prime }+\alpha _{1}t_{2\left( n-1\right)
}^{\prime }, \\ 
k_{2\left( n-1\right) } & = & 2\left( n-1\right) \alpha _{4}t_{2\left(
n-1\right) }+\left( n-1\right) \alpha _{3}t_{\left( n-1\right) }+\alpha
_{2}t_{2\left( n-1\right) }^{\prime }+\alpha _{1}t_{\left( n-1\right)
}^{\prime }, \\ 
k_{n-1} & = & \left( n-1\right) \alpha _{4}t_{\left( n-1\right) }+\alpha
_{2}t_{\left( n-1\right) }^{\prime }+\alpha _{1}t_{0}^{\prime }, \\ 
k_{0} & = & \alpha _{2}t_{0}^{\prime }.%
\end{array}
\label{Eq15A}
\end{equation}%
Compute and rewrite $k_{i}$ given by Eq. $\left(\ref{Eq15A}\right) $ in more
clear form%
\begin{equation}
\begin{array}{lll}
k_{3\left( n-1\right) } & = & \left( \alpha _{5}\left( n-2\right) \lambda
^{3}\sin ^{n+5}R+\sum_{i=0}^{8}m_{i}\left( \phi ,\eta \right) \sin
^{i}R\right) \cos R, \\ 
k_{2\left( n-1\right) } & = & \left( \alpha _{6}\left( n-2\right)
^{2}\lambda ^{2}\sin ^{2n+2}R+\sum_{i=0}^{n+5}n_{i}\left( \phi ,\eta \right)
\sin ^{i}R\right) \cos R, \\ 
k_{n-1} & = & \left( \alpha _{7}\left( n-2\right) ^{3}\lambda \sin
^{3n-1}R+\sum_{i=0}^{2n+2}p_{i}\left( \phi ,\eta \right) \sin ^{i}R\right)
\cos R, \\ 
k_{0} & = & \left( \alpha _{8}\left( n-2\right) ^{4}\sin
^{4n-4}R+\sum_{i=0}^{3n-1}q_{i}\left( \phi ,\eta \right) \sin ^{i}R\right)
\cos R,%
\end{array}
\label{Eq16A}
\end{equation}%
Here, $m_{i},n_{i},p_{i},q_{i}$ denotes the polynomials, respectively, in $%
\phi $ and $\eta ,$ and the constants represent%
\begin{eqnarray*}
\alpha _{5}
&=&2n^{7}-44n^{6}+325n^{5}-807n^{4}-796n^{3}+6906n^{2}-10\,227n+4545, \\
\alpha _{6} &=&10n^{8}-273n^{7}+3089n^{6}-18\,843n^{5}+67\,223\allowbreak
n^{4}-140\,907n^{3} \\
&&+162\,003n^{2}-81\,929n+5851, \\
\alpha _{7}
&=&3n^{9}-104n^{8}+1549n^{7}-13\,028n^{6}+68\,261n^{5}-230\,910n^{4} \\
&&+502\,291n^{3}-670\,392n^{2}+486\,104n-137\,246, \\
\alpha _{8} &=&-3\left( n-7\right) \left( n-3\right) \left( n-1\right)
\left( n^{4}-24n^{3}+194n^{2}-624n+709\right) .
\end{eqnarray*}%
Eliminating $f^{3\left( n-1\right) }$ by Eqs. $\left( \ref{Eq14}\right) $
and $\left( \ref{Eq15}\right) $, the following appears
\begin{equation}
h_{2\left( n-1\right) }f^{2\left( n-1\right) }+h_{n-1}f^{n-1}+h_{0}=0,
\label{Eq17}
\end{equation}%
where 
\begin{eqnarray*}
h_{2\left( n-1\right) } &=&t_{2\left( n-1\right) }k_{3\left( n-1\right)
}-k_{2\left( n-1\right) }t_{3\left( n-1\right) }, \\
h_{n-1} &=&t_{n-1}k_{3\left( n-1\right) }-k_{\left( n-1\right) }t_{3\left(
n-1\right) }, \\
h_{0} &=&t_{0}k_{3\left( n-1\right) }-k_{0}t_{3\left( n-1\right) }.
\end{eqnarray*}%
Using $t_{i}$ by Eq. $\left( \ref{14a}\right) $ and $k_{i}$ by Eq. $\left( %
\ref{Eq16A}\right) $ for $i=0$, $n-1$, $2\left( n-1\right) $, $3\left(
n-1\right) $,  compute $h_{j}=h_{j}(R)$ for $j=0$, $n-1$, $2\left(
n-1\right) $:%
\begin{equation}
\begin{array}{lll}
h_{2\left( n-1\right) } & = & \left( \varsigma _{1}\left( n-2\right)
^{2}\lambda ^{5}\sin ^{2n+8}R+\sum_{j=0}^{n+11}p_{2\left( n-1\right)
j}\left( \phi ,\eta \right) \sin ^{j}R\right) \cos R, \\ 
h_{n-1} & = & \left( \varsigma _{2}\left( n-2\right) ^{3}\lambda ^{4}\sin
^{3n+5}R+\sum_{j=0}^{2n+8}p_{\left( n-1\right) j}\left( \phi ,\eta \right)
\sin ^{j}R\right) \cos R, \\ 
h_{0} & = & \left( \varsigma _{3}\left( n-2\right) ^{4}\lambda ^{3}\sin
^{4n+2}R+\sum_{j=0}^{3n+5}p_{0j}\left( \phi ,\eta \right) \sin ^{j}R\right)
\cos R,%
\end{array}
\label{Eq17A}
\end{equation}%
where $p_{ij}\left( \phi ,\eta \right) $ $($for $j=0$, $n-1$, $2\left(
n-1\right) )$ denote polynomials in $\phi $ and $\eta ,$ and 
\begin{eqnarray*}
\varsigma _{1} &=&\alpha _{5}\left( 2n^{4}-29n^{3}+129n^{2}-219n+105\right)
+\alpha _{6}\left( n+1\right) ^{2}, \\
\varsigma _{2} &=&\alpha _{5}\left(
3n^{5}-56n^{4}+398n^{3}-1380n^{2}+2367n-1604\right) +\alpha _{7}\left(
n+1\right) ^{2}, \\
\varsigma _{3} &=&\alpha _{5}\left( n^{4}-24n^{3}+194n^{2}-624n+709\right)
+\alpha _{8}\left( n+1\right) ^{2}
\end{eqnarray*}%
describe some constants. Replacing $f^{2\left( n-1\right) }$ by Eq. $\left( \ref%
{Eq14}\right) $ with $f^{2\left( n-1\right) }=-\frac{h_{n-1}}{h_{2\left(
n-1\right) }}f^{n-1}-\frac{h_{0}}{h_{2\left( n-1\right) }}$ given by Eq. $%
\left( \ref{Eq17}\right)$, then the following occurs
\begin{equation*}
\left( -\frac{t_{2\left( n-1\right) }h_{n-1}}{h_{2\left( n-1\right) }}%
-t_{n-1}\right) f^{n-1}+t_{0}-\frac{t_{2\left( n-1\right) }h_{0}}{h_{2\left(
n-1\right) }}=0.
\end{equation*}%
Using more transparent notation, 
\begin{equation}
d_{n-1}f^{n-1}+d_{0}=0,  \label{Eq18}
\end{equation}%
where 
\begin{eqnarray*}
d_{n-1} &=&\left( n-2\right) ^{4}\left( \varsigma _{1}\varsigma
_{4}+\varsigma _{2}\varsigma _{5}\right) \lambda ^{6}\sin ^{4n+8}R, \\
d_{0} &=&\left( n-2\right) ^{5}\left( \varsigma _{1}\varsigma _{6}+\varsigma
_{3}\varsigma _{5}\right) \lambda ^{5}\sin ^{5n+5}R,
\end{eqnarray*}%
and%
\begin{eqnarray*}
\varsigma _{4} &=&3n^{5}-56n^{4}+398n^{3}-1380n^{2}+2367n-1604, \\
\varsigma _{5} &=&-2n^{4}+29n^{3}-129n^{2}+219n-105, \\
\varsigma _{6} &=&n^{4}-24n^{3}+194n^{2}-624n+709.
\end{eqnarray*}%
Replacing $f^{n-1}$ denoted by Eq. $\left( \ref{Eq14}\right) $ with $%
f^{n-1}=-d_{0}/d_{n-1}$ determined by Eq. $\left( \ref{Eq18}\right) ,$ the following appears
\begin{equation}
t_{3\left( n-1\right) }\left( -\frac{d_{0}}{d_{n-1}}\right) ^{3}+t_{2\left(
n-1\right) }\left( -\frac{d_{0}}{d_{n-1}}\right) ^{2}+t_{n-1}\left( -\frac{%
d_{0}}{d_{n-1}}\right) +t_{0}=0.  \label{Eq19}
\end{equation}%
Clearly,
\begin{equation}
-t_{3\left( n-1\right) }\left( d_{0}\right) ^{3}+t_{2\left( n-1\right)
}\left( d_{0}\right) ^{2}d_{\left( n-1\right) }-t_{\left( n-1\right)
}d_{0}\left( d_{\left( n-1\right) }\right) ^{2}+t_{0}\left( d_{\left(
n-1\right) }\right) ^{3}=0.  \label{Eq20}
\end{equation}%
Compute each terms of Eq. $\left( \ref{Eq20}\right) $:%
\begin{equation}
\begin{array}{lll}
-t_{3\left( n-1\right) }\left( d_{0}\right) ^{3} & = & \beta _{1}\text{ }%
\lambda ^{18}\sin ^{15n+21}R+o(\sin R), \\ 
t_{2\left( n-1\right) }\left( d_{0}\right) ^{2}d_{\left( n-1\right) } & = & 
\beta _{2}\text{ }\lambda ^{18}\sin ^{15n+21}R+o(\sin R), \\ 
-t_{\left( n-1\right) }d_{0}\left( d_{\left( n-1\right) }\right) ^{2} & = & 
\beta _{3}\text{ }\lambda ^{18}\sin ^{15n+21}R+o(\sin R), \\ 
t_{0}\left( d_{2\left( n-1\right) }\right) ^{3} & = & \beta _{4}\text{ }%
\lambda ^{18}\sin ^{15n+21}R+o(\sin R),%
\end{array}
\label{Eq21}
\end{equation}%
where%
\begin{equation}
\begin{array}{lll}
\beta _{1} & = & \left( n-2\right) ^{15}\left( n+1\right) ^{2}\left(
\varsigma _{1}\varsigma _{6}\mathfrak{+}\varsigma _{3}\varsigma _{5}\right)
^{3}, \\ 
\beta _{2} & = & -\left( n-2\right) ^{15}\varsigma _{5}\left( \varsigma
_{1}\varsigma _{4}+\varsigma _{2}\varsigma _{5}\right) \left( \varsigma
_{1}\varsigma _{6}\mathfrak{+}\varsigma _{3}\varsigma _{5}\right) ^{2}, \\ 
\beta _{3} & = & -\left( n-2\right) ^{15}\varsigma _{4}\left( \varsigma
_{1}\varsigma _{4}+\varsigma _{2}\varsigma _{5}\right) ^{2}\left( \varsigma
_{1}\varsigma _{6}\mathfrak{+}\varsigma _{3}\varsigma _{5}\right) , \\ 
\beta _{4} & = & \left( n-2\right) ^{15}\varsigma _{6}\left( \varsigma
_{1}\varsigma _{4}+\varsigma _{2}\varsigma _{5}\right) ^{3},%
\end{array}%
\end{equation}%
and $o(\sin R)$ denotes the lower degree terms in $\sin R.$ Thus, the following holds
\begin{equation}
\begin{array}{l}
-t_{3\left( n-1\right) }\left( d_{0}\right) ^{3}+t_{2\left( n-1\right)
}\left( d_{0}\right) ^{2}d_{n-1}-t_{\left( n-1\right) }d_{0}\left(
d_{n-1}\right) ^{2}+t_{0}\left( d_{n-1}\right) ^{3} \\ 
=\left( \beta _{1}+\beta _{2}+\beta _{3}+\beta _{4}\right) \lambda ^{18}\sin
^{15n+21}R+o(\sin R).%
\end{array}
\label{Eq22}
\end{equation}%
Since $R(I)$ contains an interval with $\left( \ref{Eq20}\right) ,$ then $%
\lambda $ should be $0.$ Therefore, $\phi =\eta ,$ i.e., $\mathcal{A}%
=\lambda \mathcal{I}_{n}.$ It is the end of proof.
\end{proof}

Now, we prove the Main Theorem in Section 1.

\begin{proof} 
A rotational hypersurfaces family $\mathfrak{%
x}$ obtained by rotating the unit speed curve $\gamma \left( r\right)
=\left( f(r),0,\hdots,0,\varphi \left( r\right) \right) ,$ with $f(u)>0$
around axis $x_{n}$ which is defined on an interval $I$. Suppose that the
Gauss map $\mathcal{G}$ \textit{of} $\mathfrak{x}$ satisfies $\mathbb{L}%
_{n-3}\mathcal{G}=\mathcal{AG}$ with a $n\times n$ matrix $\mathcal{A}%
\mathbf{.}$ For a function $R=R(r)$ satisfying $\left( f^{\prime
}(r),\varphi ^{\prime }\left( r\right) \right) =(\cos R(r),\sin R(r)),$ let
us put $J=\left\{ r\in I\mid R^{\prime }\left( r\right) \neq 0\right\} .$

Let us examine the following two scenarios: the case when $J=\emptyset $ and
the case when $J\neq \emptyset $. Now, let's observe the outcomes.

$(a).$ Assume that $J$ is nonempty. From the proof of Lemma 2,
$\mathcal{A}=$ diag$\left( \eta ,\hdots,\eta ,\phi \right) ,$ where $%
\lambda =\phi -\eta $ for constants $\eta ,\phi .$ When $\lambda =0,$ then $%
\eta =\phi .$ Thus, it follows from $\left( \ref{Eq3}\right) $ and $\left( %
\ref{Eq4}\right) $ that $s_{n-2}=c$ is constant, and then the
Gauss--Kronecker curvature and the mean curvature satisfy $K-cH=\phi /\left(
n-1\right) .$

When $\phi \neq 0,$ $H$ and $K$ are nonzero constants. Then, $\mathfrak{x}$
is an open part of a hypersphere (see \cite{Levi}, for 3-space case). Using $%
\left( \ref{Eq7}\right) $ and $\left( \ref{Eq13}\right) $ with $\lambda =0,$
and then $R^{\prime }$ is constant, and $f(r)=c\sin R,$ where $c\in 
\mathbb{R}^{+}.$ That is, the profile curve $\gamma $ is an open part of a
half circle centered on the rotation axis of $\mathfrak{x}$. Then, $%
\mathfrak{x}$ is an open portion of a round hypersphere.

When $\phi =0,$ then $cH=K.$ Assume $K=0,$ then $H=0$. Thus, rotational
hypersurface is flat and minimal. Therefore, $\mathfrak{x}$ is an open part
of a hyperplane. When $K\neq 0,$ then $H\neq 0.$ Hence, $\mathfrak{x}$ is an
open part of a right circular hypercone, or a circular hypercylinder.

$(b).$ Assume that $J$ is empty. Then, the profile curve $\gamma $
is a straight line. Therefore, $\mathfrak{x}$ is an open part of a
hyperplane, a right circular hypercone, or a circular hypercylinder. The
converse is clear from $\left( \ref{L_(n-3) G}\right).$
\end{proof}

As a consequence, we possess the subsequent characterizations.

\begin{corollary} 
\textit{Let} $\mathfrak{x}$ \textit{be a rotational
hypersurface in }$\mathbb{E}^{n}$\textit{. Then, following are equivalent:}

$(a)$ $\mathfrak{x}$ \textit{is an open part of a round hypersphere.}

$(b)$\textit{\ The Gauss map }$\mathcal{G}$\textit{\ of }$\mathfrak{x}$%
\textit{\ satisfies }$\mathbb{L}_{n-3}\mathcal{G}=\mathcal{AG}$ \textit{for
a regular} $n\times n$ \textit{matrix} $\mathcal{A}\mathbf{.}$
\end{corollary} 

\begin{corollary} 
\textit{Let} $\mathfrak{x}$ \textit{be a rotational
hypersurface in }$\mathbb{E}^{n}$\textit{. Then, following are equivalent:}

$(a)$ $\mathfrak{x}$ \textit{is an open part of a round right circular
hypercone.}

$(b)$\textit{\ The Gauss map }$\mathcal{G}$\textit{\ of }$\mathfrak{x}$%
\textit{\ satisfies }$\mathbb{L}_{n-3}\mathcal{G}=\mathcal{AG}$ \textit{for a%
} $n\times n$ \textit{matrix} $\mathcal{A}$.
\end{corollary}

\end{document}